\magnification \magstep1 
\input epsf.sty
\input graphicx
\input  amssym.tex
\def\sqr#1#2{{\vcenter{\hrule height.#2pt              
     \hbox{\vrule width.#2pt height#1pt\kern#1pt
     \vrule width.#2pt}
     \hrule height.#2pt}}}
\def\square{\mathchoice\sqr{5.5}4\sqr{5.0}4\sqr{4.8}3\sqr{4.8}3}
\def\qed{\hskip4pt plus1fill\ $\square$\par\medbreak}

\def\cC{{\cal C}}
\def\cF{{\cal F}}
\def\cG{{\cal G}}
\def\cX{{\cal X}}
\def\cY{{\cal Y}}

\def\P{{\bf P}}


\centerline{\bf  Continuous Families of Rational Surface}

\centerline{\bf  Automorphisms with Positive Entropy}
\medskip
\centerline{Eric Bedford\footnote*{Supported in part by the NSF} and Kyounghee Kim}

\bigskip

\noindent{\bf  \S0. Introduction. }   
Cantat [C1] has shown that if a compact projective surface carries an automorphism of positive entropy, then it has a minimal model which is either a torus, K3, or rational (or a quotient of one of these).  It has seemed that rational surfaces which carry automorphisms of positive entropy are relatively rare.  Indeed, the first infinite family of such rational surfaces was found only recently (see [BK1,2] and [M]).  Here we will show, on the contrary, that positive entropy rational surface automorphisms are more ``abundant'' than the torus and K3 cases, in the sense that they are contained in families of arbitrarily high dimension. 

We define our automorphisms in terms of birational models.  We say that a birational map $f$ of ${\bf P}^2$ is an automorphism if there are a rational surface $\cX=\cX_f$, an iterated blowup map $\pi:\cX\to{\bf P}^2$, and an automorphism $f_\cX$ of $\cX$ such that $\pi\circ f_\cX=  f\circ\pi$.   We will consider birational maps of the form
$$f(x,y)=\left(y,-x+cy+ {\sum_{{\scriptstyle \ell=2\atop {\scriptstyle \ell \ {\rm even}}}}^{k-2} }{a_{\ell}\over y^{\ell}} +{1\over y^{k} } \right)\eqno(0.1)$$
where the sum is being taken only over even values of $\ell$. 

\proclaim Theorem 1.  Let $1\le j<n$ satisfy $(j,n)=1$.  There is a nonempty set $C_n\subset {\bf R}$ such that for even $k\ge2$ and for  all choices of $c\in C_n$ and $a_{\ell}\in{\bf C}$, the map $f$ in  (0.1)
is an automorphism with entropy $\log\lambda_{n,k}$, where $\lambda_{n,k}$ is the largest root of the polynomial 
 $$\chi_{n,k}(x) = 1-k \sum_{\ell=1}^{n-1}x^\ell + x^n.\eqno(0.2)$$

The set $C_n$ will be defined in \S1; when $n$ is even, it consists of all values  $2\cos(j\pi/n)$ for $(j,n)=1$.  We will show in \S5 that the family $f$ in (0.1)  varies nontrivially with the  parameters $a_{j}$.   If the $a_j$'s are real, then $f$ is an area-preserving automorphism of the real surface obtained by taking the real points of  ${\cal X}$.  Figure 1 shows an example of the case $k=4$, $n=2$.  The large disk represents the real projective space (the real points of ${\bf P}^2$) and the bounding circle is the line at infinity.  The manifold ${\cal X}$ is obtained by performing blowing 9 blowups over the points $[0:1:0]$ and $[0:0:1]$ in the line at infinity.  The map $f$ in Figure 1 has 3 real fixed points.  Two of them are saddles, and long arcs inside the unstable manifolds are shown.  Mappings of the form (0.1) are reversible: $f$ is conjugate to $f^{-1}$ via the involution $\tau(x,y)=(y,x)$.  Thus the corresponding picture of stable manifolds would be obtained by flipping the picture about the line $\{x=y\}$.  The third fixed point, in the small island on the lower left, is elliptic, and the map shows twist-map behavior there.  Several orbits showing KAM curves and island chains have been greatly magnified and are displayed off to the left hand side.

\centerline{\includegraphics[width=1.8in]{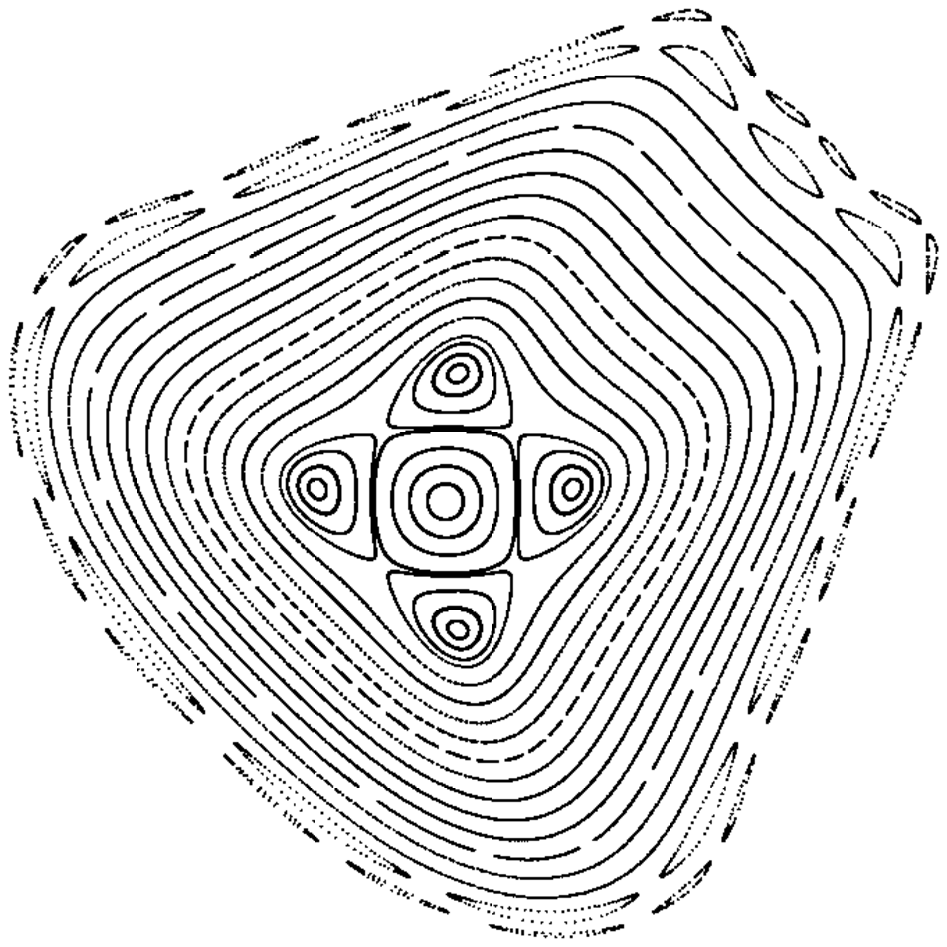} \ \includegraphics[width=2.5in]{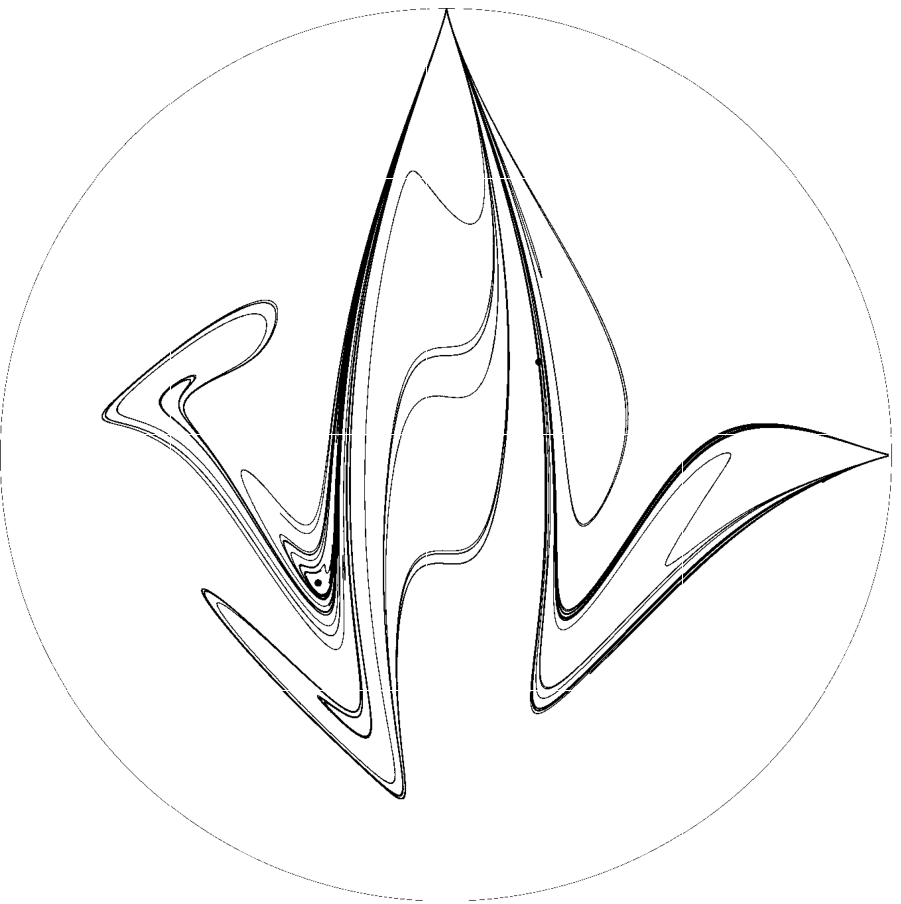}  }

\centerline{Figure 1.  $f(x,y)=(y,-x + a/y^2 + 1/y^4)$  \ \ \ $a=-2.64$}
\medskip
Each of these maps is associated with a rational surface $\cX=\cX_f$ on which the map acts as an automorphism.  We can ask to know the full automorphism group of $\cX$.  We let $Pic(\cX)$ denote the Picard group of integral divisors on $\cX$ modulo linear equivalence.   We let $\cC(\cX)$ denote the group of Cremona isometries, that is, elements of $GL(Pic(\cX))$ which are isometries with respect to the intersection product, and which preserve both the canonical class of $\cX$ and the set of effective divisors (see [D1,2] and [S]).   We discuss $Aut(\cX_f)$ via the representation $cr:Aut(\cX_f)\to\cC(\cX_f)$ given by $cr(\varphi)=\varphi_*$.  In all of our cases, this representation is  finite-to-one (Proposition 7.1); and if $a_{k-2}\ne0$ in (0.1), then $cr$ is faithful.  We define the set of effective isometries as $\cC_{ef}:=cr(Aut(\cX))$.
When $n=2$ the automorphism group of $\cX$ is maximal in the sense that every Cremona isometry is realized:
\proclaim Theorem 2.  Let $n=2$, let $f$ be as in Theorem 1, and let $\rho$ denote the reflection $(x,y)\mapsto(y,x)$.  Then $\cC_{ef}(\cX)=\cC(\cX)$ is the infinite dihedral group with generators $f_*$ and $\rho_*$.

We conclude with the observation that when $n=2$ the manifolds $\cX_a$ are generically biholomorphically inequivalent.
\proclaim Theorem 3.  Let  $n=2$ and $k\ge4$ be even; and for $a\in{\bf C}^{{k\over 2}-1}$, let $f_a$ be the map in (0.1) with corresponding manifold $\cX_a$.  There is a neighborhood $U$ of 0 in ${\bf C}^{{k\over 2}-1}$  such that if $a,\hat a\in U$, $a\ne\hat a$, and $a_{k-1}\ne0$, then $\cX_a$  is not biholomorphically equivalent to  $\cX_{\hat a}$.

${\cal C}_{ef}$ is contained in the Weyl group $W_N$, and it is interesting to know which elements of $W_N$ are actually realized by automorphisms.  McMullen [M] showed that if $w\in W_N$ has spectral radius bigger than 1 and no periodic roots, then $w$ may be realized by a rational surface automorphism $g$, i.e., $w=cr(g)$.  For the map $f$ defined by (0.1), on the other hand, $f_*$ has periodic roots.  In \S4 we factor $f_*$ as an element of two different reflection groups.  

The maps  we introduce here, as well as our methods, are  motivated by the map of Hietarinta and Viallet [HV1,2] and the subsequent study of that map by Takenawa [T1--3].  In fact, the Hietarinta-Viallet map is the map (0.1) in the case $n=3$, $k=2$.

A difference between the situation here and that of [BK1,2] and [M] is that the earlier manifolds were made by simple blowups, whereas the present ones require iterated blowups.  Some of the maps defined in [BK1,2] have invariant curves, but it seems that most of them do not.   For the maps  (0.1), the line at infinity $\Sigma_0$ is an invariant curve (the restriction to ${\Sigma_0}$ is a rotation of period $n$).  In Theorem 8.1 we show that most of the blowup fibers of the curves at infinity are curves of points where $f^{2n}$ is  tangent to the identity.

This paper is organized as follows:  \S1--3 are devoted to proving Theorem 1.   This lets us determine the invariant curves (Theorem 3.5), and we then show (Theorem 3.6) that $(f,\cX_f)$ is  minimal.   \S4 discusses reflection groups and decomposes $f_*$ into a product of reflections.  \S5 gives the nontrivial dependence of the family (0.1) on  all of the parameters.     \S6 gives properties of the Cremona isometries of $\cX$.  \S7 discusses $Aut(\cX)$ in terms of the effective Cremona isometries $\cC_{ef}(\cX)$; this is then used to prove Theorems 2 and 3.  In \S8 we show that $f^{2n}$ is tangent to the identity on the line at infinity, as well as most of the blowup fibers.

\bigskip\noindent{\bf \S1.  Construction of ${\cal X}$. } Let us write a point of ${\bf P}^2$ as $[x_0:x_1:x_2]$, and  imbed ${\bf C}^2$ into ${\bf P}^2$ via the map $(x,y)\mapsto [1:x:y]$.   We may describe the behavior of $f$ on ${\bf P}^2$ as follows.  There is a unique point of indeterminacy $e_1=[0:1:0]$ for $f$ and a unique exceptional curve $\Sigma_2=\{x_2=0\}\mapsto e_2=[0:0:1]$.  There is also a unique exceptional curve for $f^{-1}$: $ \Sigma_1=\{x_1=0\}\mapsto e_1$.  The line at infinity $\Sigma_0:=\{x_0=0\}$ is invariant.  If we write points of $\Sigma_0-e_2$  as $[x_0:x_1:x_2]=[0:1:w]$, then $f[0:1:w] = [0:1: c-1 /w]$.  Thus $f|_{\Sigma_0}$ is equivalent to the linear fractional transformation $g(w):=c- 1/w$.  

If $g$ is periodic with period $n$, then at each fixed point $w_{\rm fix}$ of $g$, we will have $g'(w_{\rm fix})=1 /w_{\rm fix}^{2}=e^{2\pi i j/n}$ for some $j$ which is relatively prime to $n$.  The equation $g(w)=w$ for a fixed point gives $w_{\rm fix}=\left(c\pm\sqrt{c^2-4}\right)/2$.  Thus set of values  $c$ for which $g$ has period $n$ is exactly
$$\{2\cos(j\pi/n):0<j<n, (j,n)=1\}.\eqno(1.2)$$
Let us use the notation $w_s=g^{s-1}(c)$ for $1\le s\le n-1$.  In other words, these are the $w$-coordinates for the forward orbit $f^se_2=[0:1:w_s]$, $1\le s\le n-1$.  
\proclaim Lemma 1.1.  For $1\le j\le n-2$, $w_jw_{n-1-j}=1$.  If $n$ is even, then $w_1\cdots w_{n-2}=1$.  If $n$ is odd, then we let $w_*(c)=w_{(n-1)/2}$ denote the midpoint of the orbit.  In this case, we have $w_1\cdots w_{n-2}=w_*$.


\noindent{\it Proof. }  We show first that $w_j w_{n-2-j}=1$ for all $1\le j\le n-2$.  To begin with, note that $g^{-1}(w)=1/(c-w)$.  
Since $w_1=c$ and $w_{n-1}=0$, we have $w_{n-2}=1/c$, so the assertion holds for $j=1$.  We now proceed by induction.  If $w_j w_{n-1-j}=1$, then $w_{j+1}=g(w_j)$ and $w_{n-1-(j+1)}=g^{-1}(w_{j}^{-1})$, and these two numbers multiply to $1$.  

The Lemma now follows if $n$ is even.  If $n$ is odd, we conclude from the first part that the product is $1$.  \qed

Let us define
$$C_n=\{c=\pm 2\cos(j\pi/n): (j,n)=1, w_1\cdots w_{n-2}=w_*=1\},\eqno(1.3)$$
where the $\pm$ notation means that we choose ``$+$'', ``$-$'', or both, corresponding to when the condition $w_*=1$ holds.
We let $\varphi(n)$ denote the number of integers $0<j<n$ which are relatively prime to $n$.  We obtain the following from Lemma 1.1:
\proclaim Lemma 1.2.  If $n$ is even, then $C_n$ coincides with the set (1.2), and thus $\#C_n=\varphi(n)$.  If $n$ is odd, then for $c$ in the set (1.2), we have $\{w_*(c),w_*(-c)\}=\{1,-1\}$.  Thus exactly one of the values $c$ or $-c$ will satisfy (1.3) and thus belong to $C_n$, and so $\#C_n={1\over 2}\varphi(n)$.

Let us set
$$q(x,y)=1+a_{k-1}y+a_{k-2}y^2+\cdots+a_1 y^{k-1}-xy^k+cy^{k+1}.\eqno(1.4)$$
We define  $b_i$, $0\le i\le 2k$ by setting
$$\eqalign{{y^k\over q(x,y)}  &  = y^k\sum_{i=0}^\infty\left(-(a_{k-1}y+\cdots +cy^{k+1})\right)^i\cr
&=\sum_{i=0}^{2k-1}b_i y^i +( x+b_{2k})y^{2k} + O(y^{2k+1}).\cr}\eqno(1.5)$$
From this it is evident that $b_0=\cdots=b_{k-1}=0$, $b_k=1$,  and $b_{k+\rho}=-a_{k-\rho}$, where $\rho>1$ is the smallest number for which $a_{k-\rho}\ne0$.  
\proclaim Lemma 1.3.  If the $a_j$'s are as in (0.1), that is, if $a_j\ne0$ only when $j$ is even, then we will have $b_j\ne0$ only when $j$ is even.

\noindent{\it Proof. }  By our hypothesis on the $a_i$'s, terms of the form  $a_{k-\rho}y^\rho$ in (1.4) can be nonzero only if $\rho$ is even.  Thus if $b_\tau y^{k+\tau}$ is a nonzero term in (1.5) with $\tau<k$, then $\tau$ must be even.  \qed

Now we construct $\cX$ by performing blowups in stages.
We begin  by making point blowups over the centers 
$\{e_2,fe_2,\dots,f^{n-1}e_2\}$. 
Let $\pi_1:\cX_1\to \P^2$ denote the resulting manifold, and let $\cF_s^1:=\pi_1^{-1}(f^se_2)$ denote the exceptional fibers.  
In a neighborhood of $e_2$ we use the local coordinate chart $(t,x)\mapsto[t:x:1]$.  For $\cF^1_0$ we will use $\pi_1:(t_1,\eta_1)\mapsto(t_1,t_1\eta_1)=(t,x)$.   Then for $\cF^1_s$, $1\le s\le n-1$, we use the coordinate chart $\pi_1:(t_1,\eta_1)\mapsto (t_1,\eta_1 t_1)=(t,y)=[t:1:y]$.

Now we continue with $2k$ more blowups over each fiber $\cF_s^{1}$, $0\le s\le n-1$.  We will proceed inductively in $j$.   We blow up the point $\xi_{j}=\beta\in\cF_s^{j}$, $1\le j \le 2k$, with $\beta$ to be specified below.  We will use the coordinate system
$$\pi_{j+1}:(\xi_{j+1},x_{j+1})\mapsto (\xi_{j+1}x_{j+1}+\beta,x_{j+1})=(\xi_{j},x_{j}).\eqno(1.7)$$
and write the fiber $\cF_s^{j+1}=\{x_{j+1}=0\}$.  The specific values of $\beta$ that we use as centers of blowup will vary with $s$ and $j$.  Over $e_2$, we have $\beta=\xi_{j}=b_{j}\in\cF_0^{j}$.  Over $fe_2$ we take $\beta=\eta_1=- b_1\in\cF_1^{1}$ at the first level, and $\beta=\xi_{j}=(-1)^{1-j}b_j\in\cF_1^{j}$ for $2\le j\le 2k$.  For $2\le s\le n-1$, we take $\beta=\xi_{j}=(w_1\cdots w_{s-1})^{j-2}b_{j}\in\cF_s^{j}$.  

\medskip

\centerline{\includegraphics[width=3.5in]{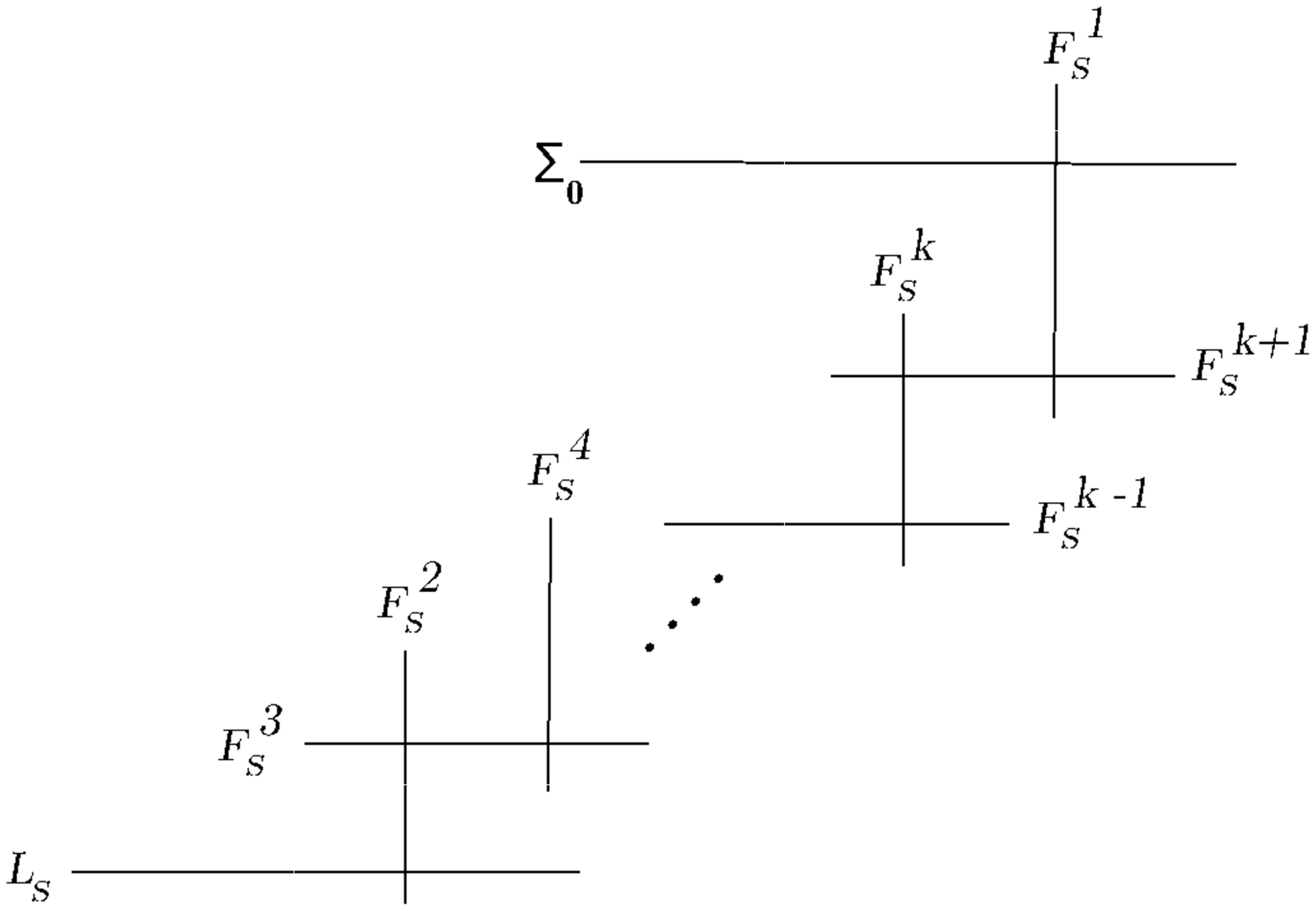}}





\centerline{Figure 2.  First sequence of blowups in construction of ${\cal X}$}

After (1.5) we saw that $b_0=\cdots=b_{k-1}=0$.  Let us interpret what this means about our space $\cX$.  We let $L_s$ denote the line in $\P^2$ connecting the origin $[1:0:0]$ to $w_s\in\Sigma_0$. Thus $L_0=\Sigma_1$ and $L_{n-1}=\Sigma_2$ (cf.\ the bottom parts of Figures 2 and 4).  The strict transform of $L_s$ inside $\cX_1$ will intersect $\cF^1_s$ at a unique point, which is the point with coordinate fiber equal to zero and which will be our center of blowup in $\cF^1_s$.  The subsequent blowup points to create $\cF^{j+1}_s$ are then taken to be $\cF^1_s\cap\cF^j_s$ for $2\le j\le k$.  At this halfway stage, we arrive at the configuration in Figure 2.  (This blowup sequence is discussed  in [BKTAM, \S2].)

\centerline{\includegraphics[height=3in]{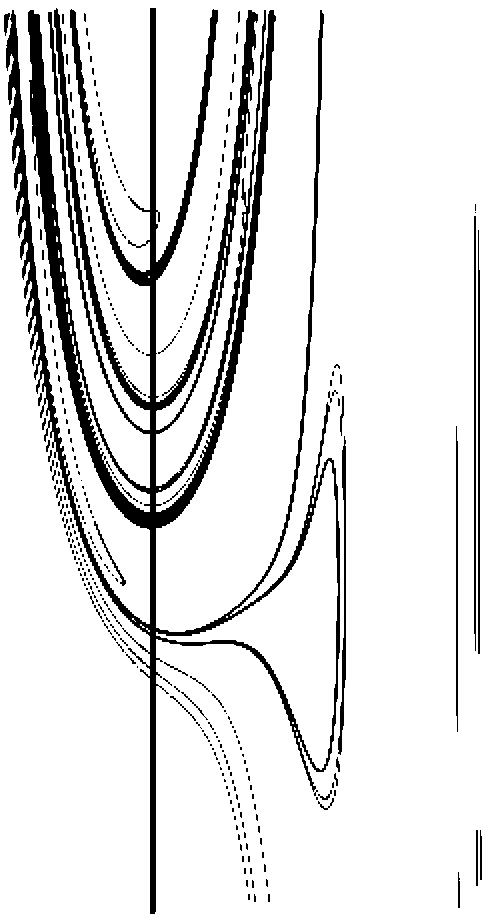}}

\centerline{  Figure 3.  Curves from Figure 1 redrawn in a neighborhood of  ${\cal F}_0^9$ (vertical axis).}
\medskip
The unstable manifold  in Figure 1 is a smooth submanifold of  ${\cal X}$.  We show in \S2 that all fibers except ${\cal F}_s^9$ are periodic.  Thus ${\cal F}_0^9$ is the only part of the fiber over $e_2$ that can intersect it.  Figure 3 shows the curve of real points of the unstable manifold from Figure 1 in $(x_9,\xi_9)$-coordinates.  The heavy, straight vertical line is the arc $\{x_9=0,2\le \xi_9\le 12\}\subset{\cal F}^9_0$.
\bigskip
\noindent{\bf \S2. Mappings of the fibers. }  Let $\cX$ be the manifold constructed in \S1, and let $f:=f_{\cX}$ denote the induced birational map of $\cX$.   We will show in Lemmas 2.1--4 that the exceptional fibers are mapped as in (2.1), with all maps being dominant. This will allow us to conclude in Proposition 2.5 that $f_\cX$ is an automorphism.
$$\eqalign{&\cF_0^1\to\cdots\to\cF_s^1\to\cF_{s+1}^1\to\cdots\to\cF_{n-1}^1\to\cF_0^1\cr
\cF_0^j\to\cdots\to\cF_s^j & \to\cF_{s+1}^j\to\cdots\to\cF_{n-1}^j\to\cF_0^{2k+2-j}\to\cdots\to\cF_{n-1}^{2k+2-j}\to\cF_0^{2k+2-j}\cr
&\Sigma_2\to\cF_0^{2k+1}\to\cdots\to\cF_{n-1}^{2k+1}\to\Sigma_1\cr}\eqno(2.1)$$

We start by seeing how $f$ maps points in the (standard) coordinate system of ${\bf C}^2$ to various coordinate neighborhoods at the fibers over $e_2$.  For instance, mapping into a coordinate system near $\cF^j_0$, we have
$$f:{\bf C}^2\ni(x,y)\mapsto (\xi_j,x_j)=(y^{k-j+1}/q(x,y),y),\ \ 2\le j\le k,$$
In the notation of (1.5), we may write the step $j=k+1$ as
$$f:(x,y)\mapsto(\xi_{k+1},x_{k+1})=(b_k+b_{k+1}y^1+\cdots + b_{2k-1}y^{k-1}+( x +b_{2k})y^{k} + O(y^{k+1}),y).\eqno(2.1)$$
Now we follow via $\pi_j^{-1}$ up to the fiber $\cF^{2k+1}_0$ and obtain:
\proclaim Lemma 2.1.  $f$ maps $\Sigma_2$  according to $f:\Sigma_2\ni[1:x:0] \mapsto\xi_{2k+1}= x+b_{2k}  \in \cF_0^{2k+1}$.

Next we determine how the fibers map forward. Let us set 
$$p(s) = a_2 s +a_3 s^2 + \cdots + a_{k-1} s^{k-2} + s^{k-1}.$$
For the rest of this section, let $f_j$ denote the mapping near the fibers $\cF_s^j$.
\proclaim Lemma 2.2.  If $a_1=0$, then $f$ maps the fibers over $e_2$ as follows:  
\vskip0pt $\cF_0^1\ni\eta_1\mapsto -\eta_1\in\cF^1_1$,
\vskip0pt $\cF^j_0\ni\xi_j\mapsto(-1)^{1-j} \xi_j\in \cF^j_1$, for $2\le j\le 2k$.

\noindent{\it Proof. }  Since $a_1=0$, we have $f[t_1:t_1\eta_1:1] = [t_1:1:c- t_1 \eta_1 + t_1^2 p(t_1)].$ Using local coordinate systems defined in \S1, we have that near $\cF^1_0$  
$$f_1 : (t_1, \eta_1) \mapsto (t_1,  \eta_1 + t_1 p(t_1)).$$ 
It follows that $f : \cF^1_0 \ni \eta_1 \mapsto - \eta_1 \in \cF^1_1.$

For $2 \le j \le k+1$, we have $b_{j-1}=0$, $\pi_j:(\xi_j, x_j) \mapsto ( \xi_j x_j, x_j) = (\xi_{j-1}, x_{j-1})$ and $f_j = \pi_j^{-1} \circ f_{j-1} \circ \pi_j$. Thus we have for $2 \le j \le k+1$,
$$f_j : ( \xi_j, x_j) \mapsto \left(  { \xi_j \over \left( -1 + \xi_j x_j^{j-2} p( \xi_j x_j^{j-1})\right)^{j-1}}\,,\, - x_j + \xi_j x_j^{j-1} p( \xi_j x_j^{j-1})\right).$$ Thus we have $f: \cF^j_0 \ni \xi_j \mapsto \xi_j/(-1)^{j-1}\in \cF^j_1$,  $2 \le j \le k+1$. 

For $k+2 \le j \le 2k$, the centers of blowup are not necessarily zero. When $ j= k+2$, the blowup center of $\cF^{k+2}_0$ is $b_{k+1}=1$. Using the previous computation the blowup center for $\cF^{k+2}_1$ is $b_{k+1}=1$. The local coordinate systems for $\cF^{k+2}_0$ and $\cF^{k+2}_1$ are $\pi_{k+2}: (\xi_{k+2},x_{k+2}) \to (\xi_{k+2} x_{k+2}+ 1, x_{k+2})$. With $f_{k+1}$ defined in the previous equation, we have 
$$\eqalign{ f_{k+2} : &( \xi_{k+2}, x_{k+2}) \mapsto\cr
& \left(  { \xi_{k+2}+O(x_{k+2}^{k-2}) \over \left( -1 +  x_{k+2}^{k-1} D(\xi_{k+2}, x_{k+2})\right)^{k+1}}\,,\,x_{k+2}\left(-1 + x_{k+2}^{k-1} D(\xi_{k+2}, x_{k+2})\right)\right)}$$
where $D(\xi,x)= (\xi x+1) p ( x^k (\xi x+1))$. Thus $f: \cF^{k+2}_0 \ni \xi_{k+2} \mapsto \xi_{k+2}$. Inductively we determine the centers of blowup and we have for $k+2\le j\le 2k$ 
$$f_{j} : ( \xi_{j}, x_{j}) \mapsto \left(  { \xi_{j}+O(x_{j}^{2k-j}) \over \left( -1 + O( x_{j}^{k})\right)^{j-1}}\,,\,x_{j}\left( -1 + O( x_{j}^{k})\right)\right).$$ 
Letting $x_j \to 0$, we have $f: \cF^j_0 \ni \xi_j \mapsto \xi_j/(-1)^{j-1}\in \cF^j_1$,  $2 \le j \le 2k$.\qed

\proclaim Lemma 2.3.  If $a_1=0$ and $1\le s\le n-2$,   $f$ maps the fibers $\cF^j_s$ to $\cF^j_{s+1}$ as follows: 
\vskip0pt $\cF_s^1\ni\eta_1\mapsto \eta_1/w_s\in\cF^1_{s+1}$,
\vskip0pt $\cF^j_s\ni\xi_j\mapsto w_s^{j-2}\xi_j  \in \cF^j_1$, for $2\le j\le 2k$,

\noindent{\it Proof. }Note that we have $f[s:1:y] = [s/y:1: -1/y + c + (s/y)^2 p (s/y)]$. With local coordinate systems $\pi_1 : (t_1, \eta_1) \mapsto [ t_1:1:t_1 \eta_1+ w_s]$ near $\cF^1_s$, the mapping near $\cF^1_s$ is given by 
$$ f_1 : (t_1, \eta_1)  \mapsto \left( { t_1 \over w_s + t_1 \eta_1}\, , \, {1 \over w_s} \eta_1+    { t_1 \over w_s + t_1 \eta_1} p\left({ t_1 \over w_s + t_1 \eta_1}\right) \right).$$
Using the same argument as in the previous Lemma, we inductively define the centers of blowup and with local coordinate systems defined in \S1, so we have for $2 \le j \le k+1$
$$f_j : (\xi_j, x_j) \mapsto  \left( {\xi_j \over ( w_s + O(x_j^j) )(1/w_s + O(x_j^{j-1}))^{j-1}}\,,\, x_j(1/w_s + O(x_j^{j-1}))\right)$$ and for $k+2 \le j \le 2k$ 
$$f_j : (\xi_j, x_j) \mapsto  \left( {\xi_j+O(x^{2k-j} )\over ( w_s + O(x_j^{k+1}) )(1/w_s + O(x_j^{k+1}))^{j-1}}\,,\, x_j(1/w_s + O(x_j^{k+1}))\right).$$
Letting $t_1 \to 0$ and $x_j \to 0$ we have the desired result. \qed

Using the same computations as in the previous two Lemmas, we also have:
\proclaim Lemma 2.4.  $f$ maps the fibers over $e_1$ as follows:
\vskip0pt $\cF^1_{n-1}\ni\eta_1\mapsto \eta_1\in\cF^1_{0}$,
\vskip0pt $\cF^{k+1-\ell}_{n-1}\ni\xi\mapsto b_{k+\ell}+1/\xi \in\cF^{k+1+\ell}_{0}$,\ \  \ \ $1\le \ell\le k-1$,
\vskip0pt $\cF^{k+1+\ell}_{n-1}\ni\xi\mapsto 1/(\xi - b_{k+\ell}) \in\cF^{k+1-\ell}_{0}$,\ \  $1\le \ell\le k-1$,
\vskip0pt $\cF^{k+1}_{n-1} \ni \xi \mapsto \xi/(\xi-1) \in \cF^{k+1}_0$
\vskip0pt $\cF^{2k+1}_{n-1}\ni\xi_{2k+1}\mapsto [1:0:\xi_{2k+1}-b_{2k}]\in\Sigma_1$.

\proclaim Proposition 2.5. If $f$ is as in Theorem 1, then $f$ is an automorphism of $\cX$.

\noindent{\it Proof.} Let us consider the complex manifold $\cX^j$ obtained by blowing up to $j$th fibers over $e_2, w_s, s=1, \dots, n-1$. Using the similar argument above Lemma 2.1, the induced birational map $f_j: \cX^j \to \cX^j$ maps $\Sigma_2$ to a fiber point $b_{j-1} \in \cF^j_0$ and the inverse map $f_j^{-1}$ maps $\Sigma_1$ to the point $ b_{j-1} \in \cF^j_{n-1}$.   Combining Lemma 1.1 and Lemmas 2.2-3, we have 
$f^{n-1}_{k+1} (b_k)   = b_k.$ 
Furthermore for all  $1\le s<(k-1)/2q$, $f^{n-1}_{k+1+2sq} (b_{k+2sq}) =  b_{k+2sq}$. Since $b_j = 0$ for all odd $j$ (see Lemma 1.3),  we have that for all $1\le j\le 2k$ $f_j^{n}$ maps $\Sigma_2 $ to the point of indeterminacy. From Lemma 2.1 and 2.4, we see that $f$ has no exceptional curve on $\cX$ and therefore $f$ is an automorphism on $\cX$.\qed

\medskip

\centerline{ \includegraphics[height=2.3in]{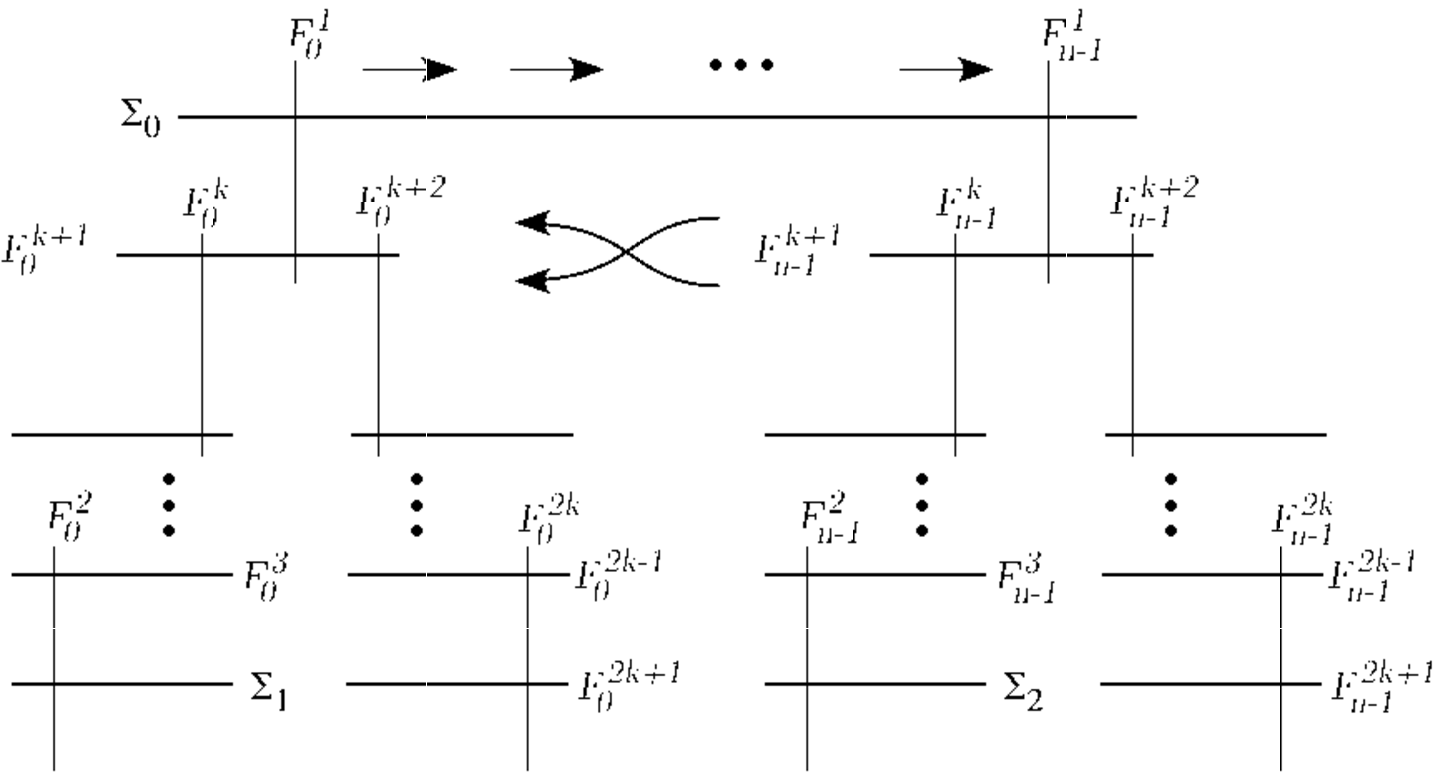}  }
\centerline{Figure 4.  Mapping the fibers. }

\medskip

Figure 4 shows graphically how the fibers are mapped, with the added information of which pairs of fibers actually intersect.  On the left, we see the fibers over $e_2$, coming off of $\Sigma_0$, starting with $\cF^1_0$.  There is a similar tree hanging off of $f^se_2$ for $1\le s\le n-1$, but the cases $1\le s\le n-2$ are not pictured.  The trees not pictured are identical, except that there is no $\Sigma_j$ connecting at the bottom.  The arrows marching to the right indicate that the arrangement hanging off of $\cF^1_s$ is mapped to the right, moving straight to the corresponding picture hanging off of $\cF^1_{s+1}$.  The twisted pair of arrows pointing to the left indicates that when we map back from $e_1$ to $e_2$, the line $\cF^{k+1}_{n-1}$ is flipped so that the fibers $\cF^{k+1\pm1}_{n-1}\to \cF^{k+1\mp1}_0$ (as well as the trees hanging below them) are interchanged.  In particular, the bottom row of the picture indicates that $\Sigma_2\to\cF^{2k+1}_0$ and $\cF^{2k+1}_{n-1}\to\Sigma_1$.  This explains the necessity for flipping because  $\cF^{2}_{n-1}$ is the only fiber that intersects $\Sigma_2$, and $\cF^{2k}_0$ is the only fiber that intersects $\cF^{2k+1}_0$.

\bigskip\noindent{\bf \S3. Behavior of $f$ on $Pic({\cal X})$. }  We will show in Corollary 3.4 that the entropy of $f$ is $\log\lambda_{n,k}$, which together with Lemma 2.5 will complete the proof of Theorem 1.  Let $S$ denote the subspace of $Pic(\cX)$ spanned by the classes of $\Sigma_0$ and $\cF^j_s$, for $1\le j\le 2k$ and $0\le s\le n-1$.  As is shown in (2.1) or in Figure 4, $S$ is invariant under $f$; indeed, $f_*$ is merely a permutation of the basis elements.  
\proclaim Proposition 3.1.  The intersection form is negative definite on $S$.

\noindent{\it Proof. }  For fixed $s$, the intersection form on the $s$-limb $\cF^1_s,\dots,\cF^{2k}_s$ (see Figure 2) is given by 
$$A_k=\pmatrix{-k-1 & 0 & 0 & 1 & 0 & 0\cr 0&-2 & 1 & &&&\cr 0&1&-2&1&&&\cr 1&&1&-2&1&&\cr 0 &&&1&-2&1\cr 0 &&&&1&-2\cr},$$
where the 1's in the first row and column are placed in the $(k+1)$st slot.  We may now write the intersection matrix $A$ on $S$ as follows.  We start with $1-n=\Sigma_0\cdot\Sigma_0$ on the upper left, and we continue down the diagonal with $n$ copies of $A_k$.  The $A_k$'s are pairwise orthothogonal, so we only need to add 1's in the first row and column at the places where $\Sigma_0\cdot\cF^1_{s}=1$.  We calculate directly that $det(A)= (1-{nk\over k+2}) [(k+2)k]^n$, so $det(A)<0$ for all the values of $k$ and $n$ that we consider.

Let $\eta_1,\dots,\eta_{1+2kn}$ denote the eigenvalues of the intersection form on $S$.  These are all nonzero since $det(A)\ne0$.  Since the intersection form has exactly one positive eigenvalue on $Pic(\cX)$, at most one of the $\eta_i$ can be positive, and the rest are negative.  However, since there is an odd number of them, and  their product is negative, we conclude that they must all be strictly negative.  \qed

Consider the action of $f_*$ on  $T:=S^\perp\subset Pic(\cX)$, the orthogonal complement with respect to the intersection product.   By Proposition 3.1, we have $S\cap T=0$ (cf. Example 5.3).  Thus $Pic(\cX)=S\oplus T$, so $dim(S)=2kn+1$ and  $dim(T)=n$.  We let $\gamma_s$ denote the projection of $\cF^{2k+1}_s$ to $T$, and thus $\{\gamma_0,\dots,\gamma_{n-1}\}$ is a basis for $T$. 

Let $\lambda_s$ denote the projection to $S$ of the line $L_s$ which connects the origin to $[0:w_s:1]$ in $\cX$.  Thus $\lambda_0$ is the projection of $\Sigma_1$, and $\lambda_{n-1}$ is the projection of $\Sigma_2$.

\proclaim Proposition 3.2.  $\lambda_s=-\gamma_s +k\sum_{t\ne s}\gamma_t$.

\noindent{\it Proof. }  For ease of notation, we work with $L_0=\Sigma_1$.  Let us start with the observation that $\Sigma_0=\Sigma_1\in Pic(\cX)$.  Pulling back by $\pi_1$ gives $\Sigma_0+\sum_{s}\cF^1_s = \Sigma_1 + \cF_0^1\in Pic(\cX^1)$ because all of the blowup centers are contained in $\Sigma_0$, but only one of them is in $\Sigma_1$.  Of the next centers of blowup, none of them are in the strict transform of $\Sigma_0$ in $\cX^1$, but one of them is $\Sigma_1\cap\cF^1_0$.  Thus we have
$$\Sigma_0+\sum_s \left(\cF^1_s+\cF^2_s\right)=\Sigma_1+\cF^1_0+ 2\cF^2_0\in Pic(\cX^2).$$
We obtain $\cX^{j+1}$ by blowing up the intersection points $\cF^j_s\cap\cF^1_s$, so we have
$$\eqalign{\Sigma_0&+\sum_s \left(\cF^1_s+\cF^2_s+2\cF^3_s+\cdots+j\cF^{j+1}_s\right)\cr
&=\Sigma_1+\cF^1_0+ 2\cF^2_0+\cdots+(j+1)\cF^{j+1}_0\in Pic(\cX^{j+1}).\cr}$$
We continue this way until we reach $\cF^{k+1}$, and thereafter we blow up free points.  This means that the coefficients stop increasing after we reach $j=k+1$, and we have
$$\Sigma_0+\sum_s\left(\cF^1_s+ \dots + k\cF^{2k+1}_s\right) = \Sigma_1+\cF_0^1+\cdots + (k+1)\cF^{2k+1}_0\in Pic(\cX).$$
This expression is a sum involving $\Sigma_1$, $\cF^{2k+1}_s$, and basis elements of $S$.  Thus, if we project everything to $T=S^\perp$, the basis elements disappear, and $\cF^{2k+1}_s$ is transformed to $\gamma_s$, from which we obtain our formula for $\lambda_0$.
\qed

By Lemma 2.4 and Proposition 3.2 we may represent the restriction $f_*|_{T}$ as
$$\lambda_{n-1}\to\gamma_0\to \gamma_1\to\cdots\to\gamma_{n-1}\to \lambda_0=-\gamma_{0}+k \gamma_1+\cdots +k \gamma_{n-1}.\eqno(3.1)$$
The characteristic polynomial of the linear transformation (3.1) is given by (0.2).
\proclaim Proposition 3.3.  The spectral radius of $f_*$ on $Pic(\cX)$ is the same as the spectral radius of $f_*|_T$ and is given by the largest root of the polynomial (0.2).

\noindent {\it Proof. }   The spectral radius of $f_*|_T$ is given by $\lambda_{n,k}$, the largest root of (0.2).  Now let $\delta(f)$ denote the spectral radius of $f_*$ on  $Pic(\cX)$.  We see from (0.2) that there is an eigenvalue $\lambda_{n,k}>1$, so $\delta(f)>1$.  Let $t$ denote an invariant class $t \in H^{1,1}({\cX})$ which is expanded by a factor of $\delta(f)$.  Since $f^*$ just permutes basis elements of $S$, it is clear that $t\notin S$.  Thus the projection of $t$ to $T$ is nonzero.  But since $t$ generates an invariant line, we must have $t\in T$.  Thus $\lambda_{n,k}=\delta(f)$.  \qed


By Cantat [C2], the entropy of an automorphism of a complex, compact surface is given by the logarithm of the spectral radius of $f_*$.  Thus we have:
\proclaim Corollary 3.4.  The entropy of $f$ is $\log\lambda_{n,k}$.

In addition to $\Sigma_0$, certain unions of the blowup fibers are invariant: these are the cycles in the first two lines in (2.1).  Conversely, there are no other invariant curves:
\proclaim Theorem 3.5.  Let $\cX_f$ be the manifold constructed from a map of the form (0.1).  Then every invariant curve is a union of components taken from $\Sigma_0$ and the blowup fibers.

 \noindent{\it Proof.}   Suppose that $\cC$ is a curve which is invariant under $f$.  Then we have a class $\cC\in Pic(\cX)$ which is invariant under $f_*$.  Let $t$ denote the orthogonal projection of $\cC$ to $T$.  This means that $f_*t=t$.  On the other hand, 1 is not a zero of $\chi_{n,k}$, so 1 is not an eigenvalue of $f_*|_T$.  Thus $t=0$.  We conclude that $\cC\in S$.  Now we know that the basis elements of $S$ are simply permuted by $f_*$, so $\cC$ must be an union of these.  \qed
 
 We say that $(f,\cX)$ is  minimal if whenever $(g,\cY)$ is an automorphism of a smooth surface, and there is a birational morphism $\varphi : \cX\to\cY$ with $\varphi\circ f=g\circ\varphi$, then $\varphi$ is an isomorphism.
 \proclaim Theorem 3.6.  $(f,\cX_f)$ is minimal if $n>2$.  If $n=2$, then it becomes minimal after we blow down $\Sigma_0$.  
 
 \noindent{\it Proof. }  Suppose that $\varphi:\cX_f\to\cY$ is a morphism.  Consider the curve $\cC$ consisting of all the varieties in $\cX_f$ which are blown down to points under $\varphi$.  It follows that $\cC$ is invariant under $f$, so by Theorem 3.5, $\cC$ must be a union of components coming from $\Sigma_0$ and $\cF^j_s$.  If $n>2$, then the self-intersection of each of the components $\Sigma_0$ and $\cF^j_s$ is $\le-2$, so it is not possible to blow any of them down.  On the other hand, if $n=2$, then the self-intersection of  $\Sigma_0$ is $-1$, so we can blow it down.  This leaves the self-intersection of all the other fibers unchanged, except for $\cF_s^1$, which increases to $-k$.  This is  strictly less than $-1$, so nothing further can be blown down.  \qed

\medskip\noindent{\bf \S4.  Two factorizations of $f_*$ into reflections. }  A basis $\{e_0,\dots, e_N\}$ for $Pic({\cal X})$ is said to be {\it geometric} if $e_0^2=1$, $e_j^2=-1$, $1\le j\le N$, and $e_i\cdot e_j=0$ for $i\ne j$.  In our case, we have $N=n(2k+1)$.  Given a geometric basis, we define roots $a_0=e_0-e_1-e_2-e_3$ and $a_j=e_{j+1}-e_j$ for $1\le j\le N-1$.  We define the corresponding reflections $\rho_j(x)=x+ (a_j\cdot x)a_j$, and we let $W_N$ denote the group generated by $r_0,\dots,r_{N-1}$.  The subgroup generated by $r_1,\dots,r_{N-1}$ is the group of permutations of $\{e_1,\dots,e_N\}$.  

We will chose a geometric basis for $Pic({\cal X})$.  Recall that ${\cal X}$ was constructed as a sequence of point blowups: ${\cal X}={\cal X}_N\to{\cal X}_{N-1}\to\cdots\to {\cal X}_0={\bf P}^2$.   Each fiber ${\cal F}^j_s$ appears for the first time in some ${\cal X}_i$ as the exceptional blowup fiber.  Writing $\pi_i:{\cal X}\to {\cal X}_i$, we let $e_s^j:=\pi_i^{*}{\cal F}_s^j$.  We let $e_0$ be the class of a general line, pulled back to ${\cal X}$.  For our geometric basis, we take $e_0$, together with $\{e_s^j\}$ for $0\le s\le n-1$ and $1\le j\le 2k+1$.  We will see that it is more convenient to give double indices to the basis elements with negative self-intersection.

Let $f$ be a map of the form (0.1).  We will write $f_*$ as an element of $W_N$.  For this, we use the notation $J$ for the Cremona inversion $r_0$, and we order our basis elements so that $a_0=e_0-e_{0}^1-e_{0}^{k+1}-e_0^{2k+1}$.  We will write $f_*$ as a composition of $J$'s and permutations of the elements $\{e_s^j\}$.  We start by defining  permutations 
$$\eqalign{  & \sigma=(0\  \dots\ (n-2)\ (n-1))\cr
& \tau=\big((2k+1)\ 2k\ (2k-1)\ \dots\  (k+2)   \big)\   \big((k+1)\ k\ (k-1)\ \dots\ 2\big)\cr 
&\phi= \phi_1\circ \phi_2, \ \ \phi_1 =(3\ \ k+1) (4\ \ k) \dots\ ({k/ 2}+3\ \ {k/ 2}+1) \cr
&\phi_2  = (2k+1\ \ k+3)(2k\  \  k+2) (2k\ \ k+2)\dots\ (3k/2 + 3\ \ 3k/2+1)\cr }$$  

We now define the permutation $\sigma_h$ to act ``horizontally'' on the basis elements, i.e., it moves from one basepoint to the next:
$$\sigma_h(e^j_s):=e^j_{\sigma (s)}.$$ 
The permutations $\tau_v$ and $\gamma_v$ act ``vertically'' in the $s=0$ fiber, i.e.,  they permute the basis elements in the fiber over $[0:1:0]\in{\bf P}^2$ and leave everything else fixed: 
$$\tau_v(e^j_{n-1})=e^{\tau (j)}_{n-1}, \ \tau_v(e^j_s)=e^j_s, \ s\ne 0, \ \ \ \phi_v(e^j_{n-1})=e^{\phi(j)}_{n-1}, \ \phi_v(e^j_s)=e^j_s, \ s\ne  0.$$
Figure 5 gives a diagram to show what is happening vertically on the fiber $s=0$.  The numbers ``1'', ``2'', \dots, refer to the geometric basis elements $e_0^1$, $e_0^2$, \dots; the large dots indicate that $J$ is the reflection generated by $e_0-e_0^1-e_0^{k+1}-e_0^{2k+1}$, and the permutations $\tau$ and $\phi$ are indicated.  

Thus we have a factorization of $f_*$:
$$f_*=  \phi_v\circ J\circ (\tau_v\circ J)^{k/2}\circ\sigma_h.$$
If we write the permutations $\sigma_h$, $\tau_v$ and $\phi_v$ as products of the transpositions $r_j$, $1\le j\le N-1$, then we will have a factorization of $f_*$ into reflections which are generators of $W_N$
\bigskip

\centerline{\includegraphics[height=1.1in]{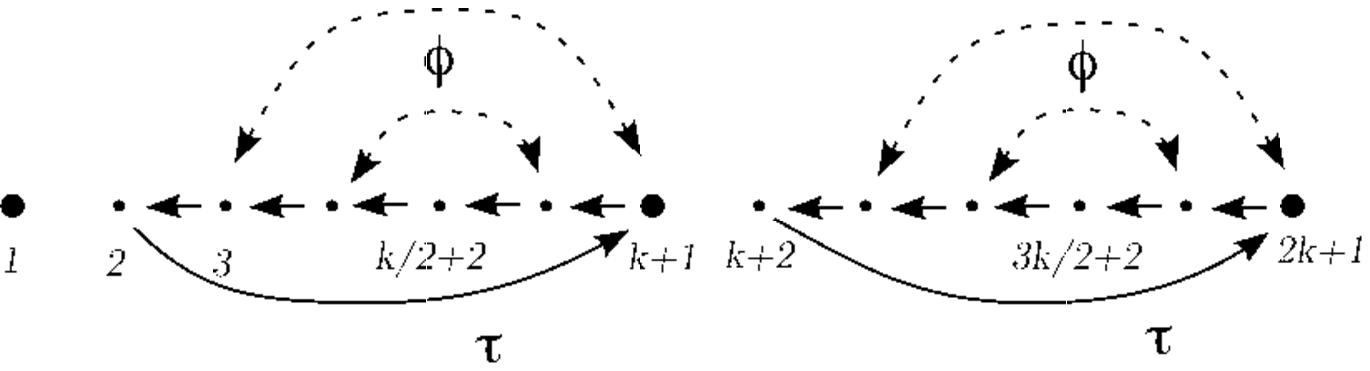} }

\centerline{Figure 5.  Permutations  $\tau$ and $\phi$.}
\medskip

Now let us see what sorts of reflections we have after we project to the subspace $T$.
 If we set $\delta:=\gamma_s\cdot\gamma_s$ and $\epsilon:=\gamma_s\cdot\gamma_t$ for $s\ne t$, then the restriction of the intersection product to $T$ is given by the matrix with $\delta$ on the diagonal and $\epsilon$ at all other places:
 $$\pmatrix{\delta & \epsilon & &\epsilon \cr
 \epsilon& \delta & \epsilon & \cr
 &\ddots&\ddots&\epsilon \cr
 \epsilon &  & \epsilon&\delta\cr}\eqno(4.1)$$
 By Proposition 3.2  the vectors $\gamma_s=(0,\dots,1,0\dots,0)$ and $\lambda_0=(-1,k,\dots,k)$ have the same self-intersection product.  It follows that up to rescaling,  we may assume that $\delta= 2-(n-2)k$ and $\epsilon=k$.

We will use the root system
$$\alpha_s:= \lambda_s-\gamma_s=(k,\dots,k,-2,k,\dots,k),\ \ 0\le s\le n-1\eqno(4.2).$$
The Cartan matrix associated with a root system is defined as the matrix $C=\left(c_{i,j}\right)$, where $c_{i,j} = 2\alpha_i\cdot\alpha_j/(\alpha_i\cdot\alpha_i)$.  By (4.1), we see that 
 $$C=\pmatrix{2 & -k  &  &-k  \cr
 -k & 2  & -k  & \cr
 &\ddots&\ddots&-k  \cr
 -k  &  & -k &2 \cr}$$
Let $\rho_s(x)=x-2{\alpha_s\cdot x\over \alpha_s\cdot\alpha_s}\alpha_s$ denote the isometric reflection generated by the root $\alpha_s$, so $\rho_s$ interchanges  $\gamma_s$ and $\lambda_s$.   In matrix form, we have 
$$\rho_{n-1} =\pmatrix{1&&&k\cr
&\ddots&&\vdots\cr
&&1&k\cr
&&&-1\cr}.$$
For $0\le s\le n-2$, let $\tau_s$ denote the permutation that transposes $\gamma_s$ and $\gamma_{s+1}$, i.e., this is the reflection defined by $\gamma_{s}-\gamma_{s+1}$.  It follows that $f_*$ is a Coxeter element of this reflection group:
$$f_*=\rho_{n-1}\circ\tau_{n-2}\circ\cdots\circ\tau_{0}.\eqno(4.3)$$  

\medskip\noindent{\bf \S5.  Nontrivial dependence on parameters. }  In this Section, we show that the family defined by (0.1) gives a $k/2-1$-dimensional family of distinct dynamical systems as we hold the parameter $c$ fixed and vary the parameters $a_\ell$.  There are $k+1$ fixed points $p_1,\dots,p_{k+1}$.    To  show that the family varies with $a_\ell$ in a nontrivial way as a family of smooth dynamical systems, it suffices to show that the eigenvalues of $Df_a$ at the point $p_j(a)$ vary with $a$.  In particular, we show that the trace of $Df$ changes nontrivially.  For this, we consider the map $a\mapsto T(a):=(\tau_1(a),\dots,\tau_{k+1}(a))$, where $\tau_j(a)$ denotes the trace of the differential $Df_a$ at $p_j(a)$.
\proclaim Proposition 5.1.  The map $a\mapsto T(a)$ has rank $k/2-1$ at the point $a=0$.

\noindent{\it Proof. }   The fixed points have the form $p_s=(\zeta_s,\zeta_s)$, where $\zeta_j$ is a root of the equation
$$\zeta=(c-1)\zeta+\sum_{{\scriptstyle j=2\atop \scriptstyle j{\rm\ even } } }^{k-2}{a_j\over \zeta^{j}}+{1\over\zeta^{k}}.\eqno(5.1)$$
When $a=0$, the fixed points all satisfy $\zeta^{k+1}=(-c+2)^{-1}$.
If we differentiate (5.1) with respect to $a_\ell$, then at $a=0$ we have $(-c+2+k/\zeta^{k+1}){\partial\zeta\over\partial a_\ell}=1/\zeta^{\ell}$, which gives $\left.{\partial\zeta\over\partial a_\ell}\right|_{a=0}=((-c+2)(k+1)\zeta^{\ell})^{-1}$.

The trace of $Df(x,y)$ is given by 
$$\tau=c-\sum_{{\scriptstyle j=2\atop \scriptstyle j{\rm\ even } } }^{k-2}{ja_j\over y^{j+1}}-{k\over y^{k+1}}.$$
If we take $y=\zeta_a$, we find that
$$\eqalign {\left.{\partial\tau(\zeta_a)\over\partial a_\ell}\right|_{a=0} &  = -{\ell\over y^{\ell+1}} +{k(k+1)\over y^{k+2}} {\partial\zeta_a\over \partial a_\ell}\cr
= -{\ell\over y^{\ell+1}} &+{k\over -c+2}{1\over \zeta^{k+1} \zeta^{\ell+1} } ={-\ell\over y^{\ell+1} }+{k\over y \zeta^{\ell}} = {k-\ell \over \zeta^{\ell+1}}\cr}$$
If we let $\zeta_j=y$ range over $k/2-1$ distinct choices of roots $(-c+2)^{{-1\over k+1}}$, then this matrix essentially is an $(k/2-1)\times(k/2-1)$ Vandermondian, so we see that it has rank $k/2-1$. \qed 

\proclaim Theorem 5.2.  Let $f_a$ be a map of the form (0.1).  There is a neighborhood $U$ of 0  in ${\bf C}^{k/2-1}$ such that if $a',a''\in U$, then $f_{a'}$ is not smoothly conjugate to $f_{a''}$.

\noindent{\it Proof. }   By Proposition 5.1, the map ${\bf C}^{k/2-1}\ni a\mapsto T(a)$ is locally injective at $a=0$.  Further, for $a=0$, the fixed points $p_s$, $1\le s\le k+1$, and thus the values $\tau_s(0)$ are distinct.  Thus the set-valued map ${\bf C}^{k/2-1}\ni a\mapsto \{\tau_1(a),\dots,\tau_{k+1}(a)\}$ is locally injective at $a=0$.  Thus if $U\ni 0$ is small, and $a',a''\in U$, $a'\ne a''$, the sets of multipliers at the fixed points are not the same, so the maps $f_{a'}$ and $f_{a''}$ cannot be smoothly conjugate.  \qed

\bigskip\noindent{\bf \S6.  Cremona Isometries. }  We will say that an element of $GL(Pic(\cX))$ is a Cremona isometry if it preserves the intersection product, if it preserves the canonical class $K_\cX$, and if preserves the set of effective divisors.  We denote the Cremona isometries by $\cC(\cX)$.   We only discuss $\cC(\cX)$ here, but we note that similar results hold for $\cC(\cY)$.

Let $\Omega$ denote the 2-form on $\cX$ which is induced from $dx\wedge dy$ on ${\bf C}^2$.  We see in the $(t,x)=[t:x:1]$-coordinates that $\Omega$ has a pole of order 3 at $\Sigma_0=\{t=0\}$.  Further, pulling back by the various coordinate maps $\pi_j$, we see that $\Omega$ has a pole at $\cF^j_s$ corresponding to the multiplicities in (6.1):
$$\eqalign{-K_{\cX}= 3\Sigma_0+\sum_{s}& \left(2 \cF^1_s+ \cF^2_s+2\cF^3_s+3\cF^4_s+\cdots+ k\cF^{k+1}_s+\right. \cr
&\left. + (k-1)\cF^{k+2}_s+(k-2)\cF^{k+3}_s+\cdots+\cF^{2k}_s\right). \cr}\eqno(6.1)$$
As in \S3, we let $S$ denote the span of $\Sigma_0$ and $\cF^j_s$, $0\le s\le n-1$, $1\le j\le 2k$ in $Pic(\cX)$.  Thus $-K_\cX\in S$.
\proclaim Proposition  6.1.  Equation (6.1) is the unique representation of the class of $-K_\cX\in Pic(\cX)$ in terms of prime divisors.

\noindent{\it Proof.  }  By Proposition 3.1 the intersection product is negative definite on $S$.  Thus the result follows from F. Sakai [S].
\qed

\medskip

\centerline{ \includegraphics[height=1.3in]{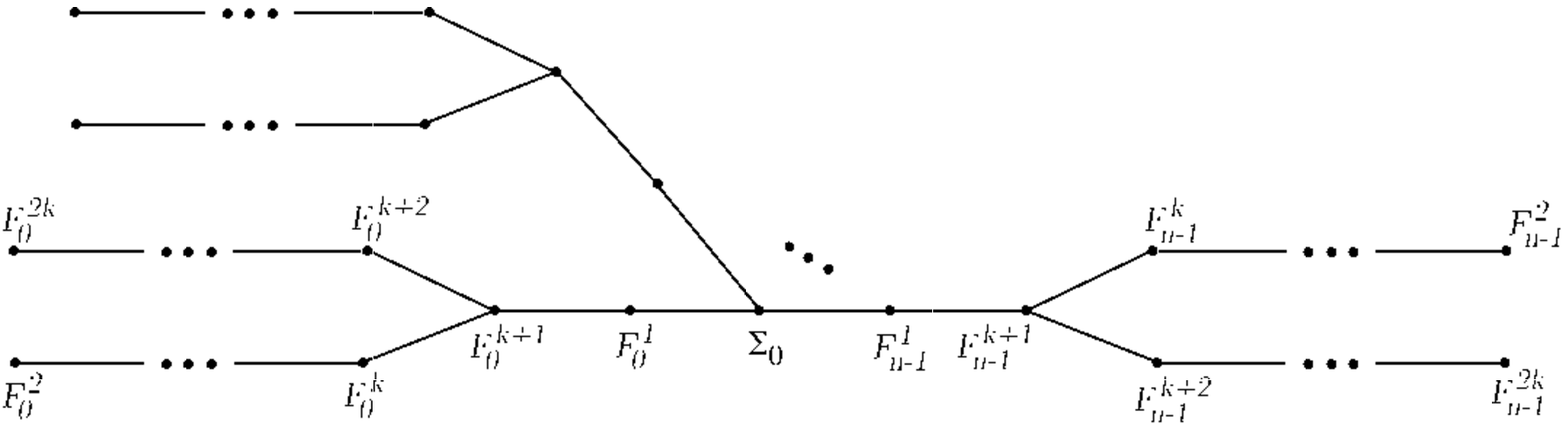}  }
\centerline{Figure 6.   Invariant graph $G_\cX$. }
\bigskip

The fibers that comprise the canonical divisor are the vertices in Figure 6.  It is dual to a portion of Figure 4: the blowup fibers are the vertices in Figure 4, and an edge indicates an intersection between the fibers.  From Proposition  6.1, we have the following:
\proclaim Proposition  6.2.   Any Cremona isometry $\cC(\cX)$ preserves the graph $G_\cX$ in Figure 6.  Thus $S$ is an invariant subspace, and the restriction $\cC(\cX)$ to $S$ is a subgroup of the permutations of the basis elements of $S$.

\noindent{\it Proof. }  A Cremona isometry preserves the anti-canonical class, so by the uniqueness of the representation (6.1),  a Cremona isometry must permute the basis elements of $S$.  Since the vertices of the graph $\cG_\cX$ represent the places where the  intersection product is  $+1$, we see that a Cremona isometry must preserve the graph.  \qed

For each $s=0,\dots, n-1$ let us set 
$$v_s := \cF_s^1+ \sum^k_{i=2} (i-1) \cF_s^i + k \sum_{i=k+1}^{2k} \cF_s^i,\ \  \ \ \ u_s := \sum^{2k}_{i=2} (i-1) \cF_s^i,$$ 
$$ \varpi_s = -k \left( \Sigma_0 + \sum_{i \ne s} v_i\right) + u_s\ \ \ \  {\rm and}\ \ \  \ \varrho_s = \varpi_s-k^2 \sum_{i\ne s} \cF_i ^{2k+1}+ 2k \,\cF_s^{2k+1}.$$
By checking through the generators of S, we see that $\varrho_s \in T$ and $\varpi_s \in S$ for $s=0,\dots , n-1$. Furthermore for $s=0, \dots, n-1$ we have  $$ k^2({k\over 2}+1) ({k\over 2}+1-n) \cF_s^{2k+1}  = \left[ ({k\over 2}+2-n) \varrho_s + \sum_{j\ne s} \varrho_j \right] - \left[ ({k\over 2}+2-n) \varpi_s + \sum_{j\ne s} \varpi_j \right] .$$
The first term in square brackets is an element of $T$, and the second term is in $S$. It follows that  $$\eqalign{\gamma_s =& -{4( k (n-3)+2(n-2) )\over k (k+2)(k-2n+2)}\cF_s^{2k+1} -{2( 4-k^2)\over k (k+2)(k-2n+2)} \sum_{j \ne s}\cF_j^{2k+1} \cr  &+{2( k-(n-2) (k^2+2k-1) )\over k^2 (k+2)(k-2n+2)}\cF_s^{2k} +{2( 4k-2-k^3)\over k^2 (k+2)(k-2n+2)} \sum_{j \ne s}\cF_j^{2k} +\cdots , \cr } \eqno{(6.2)}$$ 
where the $\cdots$ indicates a linear combination of basis elements ${\cal F}^j_s$ with $j<2k$.
The following justifies restricting Cremona isometries to $T$:
\proclaim Proposition  6.3.  A Cremona isometry of $Pic(\cX)$  is uniquely determined by its restriction to $T$.

\noindent{\it Proof. }  Suppose that $\psi\in\cC(\cX)$, and $\psi|_T$ is the identity.   By Proposition 6.2,  $\psi$ induces an automorphism of the graph $\cG_\cX$, so $\psi(\cF^{2k}_0)=\cF^{2k}_s$ for some $s$.  Using the expression in $(6.2)$ we have 
$$\cF^{2k}_0\cdot\gamma_0={2((n-4) k^2+(2n-3) k+n-2)\over k^2 (k+2)(k-2n+2)},\ \ \ \  {\rm and}$$
$$\cF^{2k}_s\cdot\gamma_0={-4(k^3-4k+1)\over k^2 (k+2)(k-2n+2)}\ \ \ \ \ \ \ \ \ \ \ \  {\rm for\ } s \ne 0.$$
Since $\psi$ is an isometry, 
$$\cF^{2k}_0\cdot\gamma_0=\psi(\cF_0^{2k})\cdot\psi\gamma_0=\psi(\cF_0^{2k})\cdot\gamma_0.$$  
It follows that $\psi(\cF^{2k}_0)=\cF^{2k}_0$.  And similarly,  $\psi(\cF^{2k}_s)=\cF^{2k}_s$ for all $s$.   Continuing down the levels, we find $\psi(\cF^j_s)=\cF^j_s$ for all $1\le j\le 2k$, and $\psi$ fixes $\Sigma_0$. Thus $\psi$ is the identity.
\qed

\noindent{\bf \S7.  The Automorphism Group of $\cX$.  }  
In this section we will describe $Aut(\cX)$ in terms of the representation $cr:Aut(\cX)\to\cC(\cX)$.

\proclaim Proposition  7.1.  Let $k$, $q$, and $f$ be as in Theorem 1, and let $\cX=\cX_f$ be the surface constructed in \S2.  If $\varphi\in Aut(\cX)$ induces the identity map on $Pic(\cX)$, then $\varphi(x,y)=(\alpha x,\beta y)$, where $\alpha^{k}\beta=\alpha\beta^{k}=1$.  In particular, $cr$ is at most $(k^2-1)$-to-one.  If   $a_{k-2}\ne0$, then $cr$ is faithful.

\noindent{\it Proof. }  If $\varphi^*$ is the identity on $Pic(\cX)$, then $\varphi$ descends to an automorphism of ${\bf P}^2$ which must fix $e_1$ and $e_2$.  Further, $\Sigma_0$, the line connecting $e_1$ and $e_2$, must be invariant.  Thus, in the coordinates $[1:x:y]$, we must have $\varphi(x,y)=(\alpha_1 x+\alpha_0,\beta_1y+\beta_0)$.

Now let us look at the fiber $\cF^1_0$; the fiber point corresponding to $\{x=0\}\cap\cF^1_0$ is a center of blowup, so it must be fixed. Thus we must have $\alpha_0=0$.  A similar argument at $e_1$ gives $\beta_0=0$.

Now in the coordinate system $(s,x)=[s:x:1]$ at $e_2$, we have $\varphi:(s,x)\mapsto(s/\beta_1,\alpha_1x/\beta_1)$.  As we pass to the various blowup coordinates, we find
$$\eqalign{\varphi:\ &(s_1,\eta_1)\mapsto (s_1/\beta_1,\alpha_1\eta_1)\cr
&(\xi_2,x_2)\mapsto (\xi_2/(\alpha_1\beta_1),\alpha_1 x_2)\cr
&\dots\cr
&(\xi_{k+1},x_{k+1})\mapsto(\xi_{k+1}/(\alpha_1^{k}\beta_1),\alpha_1 x_{k+1})\cr}$$
Note that the point $\xi_{k+1}=b_k=1$ is a center of blowup, so it must be preserved, which gives us $\alpha_1^{k}\beta_1=1$.  A similar argument at $e_1$ yields $\alpha_1\beta_1^{k}=1$.

If we substitute one of these equations into the other we find that $\alpha_1$, for instance, is a $(k^2-1)$-th root of unity.  Thus there are at most $k^2-1$ pairs $(\alpha_1,\beta_1)$.

If, in addition, we have $a_{k-2}\ne0$, then $b_{k+2}=-a_{k-2}\ne0$.  This, too, is a center of blowup, so by the same argument we have $\alpha_1^{k+2}\beta_1=1$.  Combined with the earlier equation, this gives $\alpha_1^{2}=1$.  Thus  $\alpha_1^{k}\beta_1=1$ gives $\beta_1=1$.  Similarly, we have $\alpha_1=1$, which means that $\varphi$ is the identity in this case.
\qed
%

\proclaim Proposition 7.2.  The linear map $\rho(x,y) = (y,x)$ defines an automorphism of ${\cal X}$.

\noindent{\it Proof.} We need to check that the induced map $\rho_{\cal X}$ behaves well on the various blowup fibers. If we follow through the arguments of Lemmas 2.2 and 2.3, we find that $\rho_{\cal X}$ maps fibers as follows:
$$\eqalign{ &{\cal F}_0^j \ni \xi_j  \leftrightarrow \xi_j \in {\cal F}_{n-1}^j, \ \ \ {\rm for\ } 1 \le j \le 2k+1\cr
&{\cal F}_s^j \ni \xi_j  \mapsto -(-w)^{2-j} \xi_j \in {\cal F}_{n-1-s}^j, \ \ \ {\rm for\ } 1 \le j \le 2k+1,\ 1 \le s \le n-2.\cr}$$  \qed
  
%
 
\noindent{\it Proof of Theorem 2. }  By Proposition 6.3 it suffices to consider the restriction of a Cremona isometry to $T$.  Since $n=2$, $T$ is 2-dimensional, and the intersection product (3.2) has the form $\pmatrix{2&k\cr k&2}$.  The null space $\{v\cdot v=0\}$ is generated by the vectors $v_j=(1,\lambda_j)$, $j=1,2$, with $\lambda_{{1\atop 2}}=(1,-(k\pm\sqrt{k^2-4})/2)$.  An isometry $\psi\in\cC(\cX)$ must preserve the null space, so applying the reflection $\rho_*$ if necessary, we may assume that $v_j$ is an eigenvector: $\psi v_j=s_jv_j$.  Since $\psi$ must preserve the canonical class, we have $s_j>0$, and since it preserves the lattice, we have $s_1s_2=1$.  Now we may diagonalize
$$\psi=P\pmatrix{s&0\cr 0&s^{-1}} P^{-1}$$
where $P=\pmatrix{1&1\cr \lambda_1& \lambda_2}$.   The upper right hand entry of $\psi$ is 
$$\psi_{1,2}=(-s+s^{-1})(\lambda_2-\lambda_1)^{-1}=(-s+s^{-1})/\sqrt{k^2-4}.$$
The set of all $s$ is a multiplicative cyclic group, and without loss of generality, we may suppose that $s>1$ is minimal, and thus a generator.  Now $\psi_{1,2}$ must be an integer, and the minimal value for which this can happen occurs for $\psi_{1,2}=\pm1$, in which case we have $s=(\pm \sqrt{k^2-4}\pm k)/2$.  Since $s>1$, we have $s=(k+\sqrt{k^2-4})/2$.  On the other hand, in this case we have $\psi=f_*$, which shows that  $f_*$ and $\rho_*$ generate $\cC(\cX)=\cC_{ef}(\cX)$.
\qed

\noindent{\it Proof of Theorem 3. }  Suppose that $h:\cX_a\to\cX_{\hat a}$ is a biholomorphism.  Then $g:=h^{-1}\circ f_{\hat a}\circ h\in Aut(\cX_a)$.  Recall that $T_a$ (as well as $T_{\hat a}$) has dimension 2.  The null vectors $\{v\in T_a:v\cdot v=0\}$ are eigenvectors for $f_{a*}$ (and similarly for  $f_{\hat a *}$).  It follows that $g_*$ has the same eigenvectors.   By Theorem 2, we know that $g_*$ is in the dihedral group generated by $f_{a*}$ and $\tau_*$.   Since $g_*$ and $f_{a*}$ have the same eigenvectors, we must have $g_*=f^j_{a*}$.  Since they have the same spectral radius, we must have $f_{a*}=g_*$.  By Proposition 7.1, we must have $f_a=g$.  On the other hand, this means that $f_a$ is conjugate to $f_{\hat a}$, which is not possible by Theorem 4.2.  \qed

\bigskip
\noindent{\bf \S8.  Parabolic Points. } The invariant curve $\Sigma_0$ as well as most of the blowup fibers  give  curves of parabolic points.  To see this, let us rewrite the map $f$ in $(0.1)$ near the line at infinity using the identification $(t,x) \leftrightarrow [t:x:1] \in \P^2$ $$ f(t,x) = \left( {t \over -x+c + O(t^3)}\, , \, {1 \over -x+c + O(t^3)}\right).$$ In this affine coordinate system, we have $\Sigma_0 = \{ t=0\}$ and the orbit of  each point $(0,x)$ in $\Sigma_0 \setminus \{w_s, 0 \le s \le n-1\}$ is given by $$f: (0,x) \to (0,h(x)) \to (0,h^2(x)) \to \cdots \to (0, h^{n-1}(x)) \to (0,x)$$ where $h(x) = 1/(c-x)$. Since $\Sigma_0$ is a line of fixed points, it follows that $$ Df^{n}|_{(0,x)} = \left[ \matrix{ x \cdot h(x)  \cdots h^{n-1}(x)  & 0 \cr 0 &  (x \cdot h(x)  \cdots h^{n-1}(x) )^2} \right] = \left[\matrix{\pm1&0\cr 0 & 1} \right]\eqno{(8.1)}$$
With the formulas for the mappings near blowup fibers in \S2, we have 
$$\eqalign{f: & (s, \xi ) \mapsto (s\,,\, \xi+O(s^2))  \ \ \ {\rm near\ } \cF^1_0 \cr & (s, \xi ) \mapsto ({s \over w_j+s \xi}\, ,\,  { \xi\over w_j} +O(s^2))  \ \ \ {\rm near\ } \cF^1_j,\ 1 \le j \le n-2 \cr & (s, \xi ) \mapsto (  {s \xi^k \over  - \xi^k+O(s)}\,,\,\xi)  \ \ \ {\rm near\ } \cF^1_{n-1}  .} $$
Since the differential at each point  on $\cF_s^1$ does not depend on the point, using the condition $(1.3)$ we have for each $(0,\xi) \in \cF_s^1$
$$Df^{2n}|_{(0,\xi)}  = \left[ \matrix{  (w_1 w_2 \cdots w_{n-2})^{-2} & 0 \cr 0 & (w_1 w_2 \cdots w_{n-2})^{-2} } \right] = Id\eqno{(8.2)}$$
Also we have 
$$ \eqalign{f: & (\eta, t) \mapsto ({1\ \over -1/\eta+1 + O(t^{2})}\, ,\, t) \ \ \ {\rm near\  } \cF^{k+1}_{n-1} \cr& (\eta, t) \mapsto ({1\ \over \eta+ O(t^{2})}\, ,\, t) \ \ \ {\rm near\  } \cF^{k+1\pm2i}_{n-1}, 1 \le i \le k/2-1. }$$
Furthermore we have $$\eqalign{f^n : &(\eta,0) \in \cF_{n-1}^{k+1} \mapsto ({\eta\over \eta-1},0) \in  \cF_{n-1}^{k+1} \mapsto (\eta,0) \in \cF_{n-1}^{k+1} \cr & (\eta,0) \in \cF_{n-1}^{2k+1 -2i} \mapsto ( -{1 \over \eta},0) \in \cF_{n-1}^{2k+1 +2i}\mapsto  (\eta,0) \in \cF_{n-1}^{2k+1 -2i}   .}$$
Together with formulas given in the proof of Lemma 2.2 and Lemma 2.3 we have 
$$ Df^{2n} |_{(\eta,0)} =Id , \ \ \ \ \ (\eta,0) \in \cF_s^{k+1\pm 2i} , 0 \le i\le k/2-1\eqno{(8.3)}$$
It follows that for $1 \le i \le k/2-1$ we have $f^{2n} (\eta, t ) = (\eta, t) + t^2 ( h_1(\eta, t), h_2(\eta, t))$ near $\cF_0^{k+1- 2i}$ and thus
$$f^{2n} (\eta, t) = \left( { \eta + h_1(\eta t^2, t) \over (1+ t h_2(\eta t^2, t))^2} \, ,\, t + t^2 h_2(\eta t^2, t) \right)  \ \ \ {\rm near\ \  }\cF_0^{k-1- 2i}.$$ 
Since $Df^{2n}$ is also $Id$ near $\cF_0^{k-1- 2i}$, we have $ h_1(0,0) = h_2(0,0) = 0$. Using the local coordiantes defined in \S 2, we have 
$$f^{2n} (\eta, t) = \left ( {\eta + t h_1(\eta t, t) \over 1+ t h_2(\eta t, t)}\, ,\, t+ t^2 h_2(\eta t, t) \right)\ \ \ {\rm near\ \  }\cF_0^{k- 2i}$$ 
and 
$$ Df^{2n}|_{(\eta, 0)} = \left[ \matrix{1&0\cr 0&1} \right]\  \ \ {\rm near\ \  }\cF_0^{k- 2i}, 0 \le i\le k/2-2 \eqno{(8.4)}$$ 
Combining $(8.1-4)$ we have :

\proclaim Theorem 8.1.  Let $f$ be as given in Theorem 1.  Let ${\cal P}=\Sigma_0\cup \bigcup {\cal F}_s^j$, where the union is taken over  $0\le s\le n-1$ and $j=1$ and $3\le j\le 2k-1$.  Then all points of ${\cal P}$ are fixed by $f^{2n}$, and $f^{2n}$ is tangent to the identity there.

\bigskip\noindent {\bf \S9.  Roots of unity for jacobians. }
A similar construction allows us more freedom to specify $\delta$.  In fact, it will yield mappings for which the jacobian determinant will be any root of unity.  For this we need to define the space $C(\delta,n)$, which is essentially the set of values $c=2\sqrt\delta\cos(j\pi /n)$;  we refer to [BK3] for details.  The following is similar to Theorem 1:
\proclaim Theorem 9.1.  Let $n,k,q$ be integers with $n,k\ge2$, such that  $1\le q\le k+1$ and $k+1\equiv q\ {\rm mod}\ 2q$.  Let $\epsilon\in{\bf C}$ be such that $-\epsilon^n$ is a primitive $q$th root of $-1$, and set $\delta=\epsilon^2$.  For any $c\in C(\delta,n)$, and for any choice of $a_{k-2q}, a_{k-4q},\dots\in{\bf C}$, the map
$$f(x,y)=\left (y,-\delta x+cy +\sum_{1\le s< {k-1\over 2q}} {a_{k-2sq}\over y^{k-2sq}} + {1\over y^k} \right) \eqno(9.1)$$
is an automorphism.   The entropy of this map is $\log\lambda_{n,k}$, where $\lambda_{n,k}$ is the largest root of the polynomial $\chi_{n,k}$  in (0.2).

\bigskip
\centerline{\bf References}
\medskip

\item{[BK1]}  E. Bedford and KH Kim, Periodicities in linear fractional recurrences:  Degree growth of birational surface maps, Mich.\ Math.\ J. 54 (2006), 647--670.

\item{[BK2]}  E. Bedford and KH Kim,  Dynamics of rational surface automorphisms: Linear fractional recurrences.  J. of Geometric Analysis, to appear.   arXiv:math/0611297

\item{[BK3]} E. Bedford and KH Kim, Dynamics of rational surface automorphisms: Linearization globalization.  

\item{[BKTAM]} E. Bedford, KH Kim, T. Truong, N. Abarenkova, J-M Maillard,  Degree complexity of a family of birational maps.  Mathematical Physics, Analysis and Geometry, to appear.   arXiv:0711.1186 


\item{[C]} S. Cantat, Dynamique des automorphismes des surfaces projectives complexes, C.R. Acad.\ Sci.\ Paris S\'er.\ I Math 328 (1999), 901--906.



\item{[DF]}  J. Diller and C. Favre,  Dynamics of bimeromorphic maps of surfaces, Amer. J.  Math., 123 (2001), 1135--1169.

\item{[D1]}  I. Dolgachev,  Weyl groups and Cremona transformations, Proc.\ Symp.\ Pure Math. Vol. 40 (1983) Part 1, 283--294.

\item{[D2]} I. Dolgachev, Reflection groups in algebraic geometry, Bull.\ AMS, 45 (2007), 1--60.




\item{[HV1]} J. Hietarinta and C. Viallet, Singularity confinement and degree growth. SIDE III---symmetries and integrability of difference equations (Sabaudia, 1998), 209--216, CRM Proc. Lecture Notes, 25, Amer. Math. Soc., Providence, RI, 2000.

\item{[HV2]} J. Hietarinta and C. Viallet,  Discrete Painlev\'e I and singularity confinement in projective space. Integrability and chaos in discrete systems (Brussels, 1997). Chaos Solitons Fractals 11 (2000), no. 1-3, 29--32.




\item {[M]}  C. McMullen,  Dynamics on blowups of the projective plane. Publ. Math. I.H.E.S. No. 105 (2007), 49--89.


\item{[Sa]}  F. Sakai,  Anticanonical models of rational surfaces, Math.\ Ann.\ 269, 389--410 (1984).


\item{[T1]}  T. Takenawa, A geometric approach to singularity confinement and algebraic entropy, J. Phys.\ A: Math.\ Gen.\ 34 (2001) L95--L102.

\item{[T2]} T. Takenawa, Discrete dynamical systems associated with root systems of indefinite type, Commun.\ Math.\  Phys., 224, 657--681 (2001).

\item{[T3]} T. Takenawa, Algebraic entropy and the space of initial values for discrete dynamical systems, J. Phys.\ A: Math.\ Gen.\ 34 (2001) 10533--10545.

\bigskip

\rightline{E. Bedford:  bedford@indiana.edu}

\rightline {Department of Mathematics}

\rightline{Indiana University}

\rightline{Bloomington, IN 47405 USA}

\medskip\rightline{K. Kim: kim@math.fsu.edu}

\rightline {Department of Mathematics}

\rightline{Florida State University}

\rightline{Tallahassee, FL 32306 USA}

\bye